\documentclass[reqno]{amsart}

\usepackage{amscd}

\newtheorem{thm}{Theorem}

\newtheorem{dfn}{Definition}

\newtheorem{prp}{Proposition}

\newtheorem{cor}{Corollary}

\newtheorem{lma}{Lemma}

\newtheorem{rem}{Remark}
\newenvironment{rmk}{\begin{rem}\rm}{\end{rem}}

\newtheorem{exe}{Example}
\newenvironment{exm}{\begin{exe}\rm}{\end{exe}}

\newtheorem{cnj}{Question}

\newenvironment{pf}{\begin{proof}}{\end{proof}}

\newtheorem{thmK}{Theorem}

\newcommand{\C}{{\mathbb{C}}}
\newcommand{\R}{{\mathbb{R}}}
\newcommand{\Q}{{\mathbb{Q}}}
\newcommand{\Z}{{\mathbb{Z}}}

\newcommand{\Su}{{\rm {S}}}
\newcommand{\la}{\langle}
\newcommand{\ra}{\rangle}
\newcommand{\lk}{\operatorname{lk}}
\newcommand{\Vect}{\operatorname{Vect}}
\newcommand{\Imm}{\mathbf{Imm}}
\newcommand{\Emb}{\mathbf{Emb}}
\newcommand{\img}{\operatorname{im}}
\newcommand{\rank}{\operatorname{rank}}
\renewcommand{\mod}{\operatorname{mod}}
\begin{document}
\title
[Geometric formulas for Smale invariants\dots]
{Geometric formulas for Smale invariants of codimension two immersions} 
\author{Tobias Ekholm}
\thanks{Supported by a post-doctoral grant
from Stiftelsen f{\"o}r internationalisering 
av forskning och h{\"o}gre utbildning.}
\address{Stanford University, Department of mathematics, Stanford, CA
94305}
\email{tobias@math.stanford.edu} 
\author{Andr\'as Sz\H ucs}
\thanks{Supported by grants OTKA T029759 and
FKFP 0226/1999.} 
\address{Department of Analysis, ELTE, R\'ak\'oczi \'ut 5-7, Budapest
1088, Hungary} 
\email{szucs@cs.elte.hu}
\date{}
\begin{abstract}
We give three formulas expressing the Smale invariant of an
immersion $f$ of a $(4k-1)$-sphere into $(4k+1)$-space. The terms of
the formulas are geometric characteristics of any generic smooth map $g$
of any oriented $4k$-dimensional manifold,
where $g$ restricted to the boundary is an immersion regularly
homotopic to $f$ in $(6k-1)$-space. 

The formulas imply that if $f$ and $g$ are two non-regularly homotopic
immersions of a $(4k-1)$-sphere into $(4k+1)$-space then they are also
non-regularly homotopic as immersions into $(6k-1)$-space. Moreover, 
any generic homotopy in $(6k-1)$-space connecting $f$ to $g$ must have
at least $a_k(2k-1)!$ cusps, where $a_k=2$ if $k$ is odd and $a_k=1$
if $k$ is even.
\end{abstract}

\maketitle

\section{Introduction}
Whitney ~\cite{W} classified regular plane curves up to regular
homotopy: two regular curves are regularly homotopic if and only if
they have the same tangential degree. He also gave a formula for the
tangential degree of a plane curve in terms of its double points.

Smale \cite{Sm}, generalized Whitney's result to higher dimensions:
associated to each immersion $f\colon S^k\to \R^n$, $n>k$ is its
Smale invariant $\Omega(f)\in\pi_k(V_{n,k})$, the
$k^{\rm th}$ homotopy group of the Stiefel manifold of $k$-frames in
$n$-space. Two immersions are regularly homotopic if and only if they
have the same Smale invariant.

Whitney's double point formula has straightforward generalizations to
sphere immersions in double dimension. The Smale invariant of an
immersion $S^k\to\R^{2k}$ is its algebraic number of double points
($\mod 2$ if $k$ is odd). Also in dimensions right below double
($S^k\to \R^{2k-r}$, $r=1,2$) there are double point formulas for the
Smale invariant, see ~\cite{E1} and ~\cite{E2}.

In small codimension, the first case after plane curves is immersions
$S^2\to\R^3$. Smale's work shows that they are all regularly
homotopic. Regular homotopy classes of immersions $S^3\to\R^4$ form a
group isomorphic to $\Z\oplus\Z$. A description of the two
characterizing integers is given by Hughes ~\cite{Hug}. The general
codimension one case $S^n\to\R^{n+1}$ was studied by Kaiser
~\cite{Ka}. In the present paper our main concern is the codimension
two case. Certain codimension two immersions are especially
interesting due to the following:

The groups $\pi_{4k-1}(V_{4k+1,4k-1})$, enumerating regular homotopy
classes of immersions $S^{4k-1}\to\R^{4k+1}$, are infinite cyclic. A
result of Hughes and Melvin ~\cite{HM} says that in these dimensions there
exist {\it embeddings} which are not regularly
homotopic to the standard embedding. Therefore, in contrast to the
case of high codimension, the Smale invariant
can not be expressed solely through the self intersection. However,
there are still geometric formulas for Smale invariants: 
\begin{thm}\label{thm1}
Let $f\colon S^{4k-1}\to\R^{4k+1}$ be an immersion and let
$\Omega(f)\in\Z$ be its Smale invariant. 
Let $j\colon\R^{4k+1}\to\R^{6k-1}=\partial\R^{6k}_+$ denote the
inclusion and let $M^{4k}$ be any compact oriented manifold with
$\partial M^{4k}=S^{4k-1}$.  
Let $a_k=2$ if $k$ is odd and $a_k=1$ if $k$ is even.
\begin{itemize}
\item[{\rm (a)}]
If $g\colon M^{4k}\to\R^{6k-1}$ is a
generic map such that the restriction $\partial g$ of $g$ to the
boundary is regularly homotopic to $j\circ f$, then 
\begin{align}\label{grom1}
\Omega(f)&=\frac1{a_k(2k-1)!}\left(-\bar p_k[\hat
M^{4k}]+\sharp\Sigma^{1,1}(g)\right),\\ 
&=\frac1{a_k(2k-1)!}\left(-\bar p_k[\hat
M^{4k}]+e(\xi(g))\right).\label{grom2}   
\end{align} 
\item[{\rm (b)}]
If $g\colon M^{4k}\to\R^{6k}_+$ is a generic map such that
$g^{-1}(\partial\R^{6k}_+)=\partial M^{4k}$ and such that the
restriction $\partial g$ of $g$ to the boundary is a generic immersion
regularly homotopic to $j\circ f$, then 
\begin{equation}\label{Shigh}
\Omega (f) = \frac1{a_k(2k-1)!} \left(- \bar p_k[\hat M^{4k}]+
3t(g)-3l(g)+L(\partial g)\right). 
\end{equation}
\end{itemize}
\end{thm}
\noindent
We give brief explanations of the various terms in the equations above:  

All terms appearing in  ~\eqref{grom1},
~\eqref{grom2}, and ~\eqref{Shigh} are integers.  
The term $\bar p_k[\hat M^{4k}]$ denotes the $k^{\rm th}$ normal
Pontryagin class of the closed manifold $\hat M^{4k}$, obtained by
adding a disk to 
$M^{4k}$ along $\partial M^{4k}$, evaluated on its fundamental homology
class, see Section ~\ref{basnot}. 

 The term $\sharp\Sigma^{1,1}(g)$ in ~\eqref{grom1} is the algebraic
number of cusps of $g$ and the term $e(\xi(g))$ in ~\eqref{grom2} is
the Euler number of the cokernel bundle of the differential of $g$
over its singularity set, see Section ~\ref{not1}.

In Equation ~\eqref{Shigh}, the term $L(\partial g)$ measures the linking
of the double point set of $\partial g$ with the rest of its image,
see Section ~\ref{nottrip}, the term $t(g)$ is the algebraic number of
triple points of $g$, and the term $l(g)$ 
measures the linking of the singularity set of $g$ with the rest of
its image, see Section ~\ref{not2}. 

Part (a) of Theorem ~\ref{thm1} is proved in 
Section ~\ref{pfthm1a} and part (b) in Section ~\ref{pfthm1b}.   
Equation ~\eqref{grom2} was inspired by an exercise in Gromov's book
~\cite{G}, see Remark ~\ref{exfel}.  

Let $\Imm(S^n,\R^{n+k})$ denote the set (group) of
regular homotopy classes of immersions $S^n\to\R^{n+k}$. Theorem
~\ref{thm1} implies the following (the notation is the same as in
Theorem ~\ref{thm1}): 
\begin{cor}\label{corthm1}
The natural map $\Imm(S^{4k-1},\R^{4k+1})\to\Imm(S^{4k-1},\R^{6k-1})$ is
injective. Moreover, if $f,g\colon S^{4k-1}\to\R^{4k+1}$ are two
immersions with $\Omega(f)-\Omega(g)=c\in\Z$
then the algebraic number of cusps of any generic homotopy    
$F\colon S^{4k-1}\times I\to\R^{6k-1}$ connecting $j\circ f$ to
$j\circ g$ is $c\cdot a_k(2k-1)!$. In particular, any such homotopy
has at least $|c|\cdot a_k(2k-1)!$ cusp points.
\end{cor}
\noindent
Corollary ~\ref{corthm1} is proved in Section ~\ref{pfcorthm1}.

Combining the first statement of Corollary ~\ref{corthm1} with the result
of Hughes and Melvin on embeddings $S^{4k-1}\to\R^{4k+1}$ mentioned
above, one concludes that the inequality in the following theorem of
Kervaire ~\cite{K}
\begin{thmK}
If $2q>3n+1$ then every
embedding $S^n\to\R^q$ is regularly homotopic to the standard
embedding. 
\end{thmK}
\noindent
is best possible for $n=4k-1$. (This was known to Haefliger, see
~\cite{Hae}.)  

As another consequence of Theorem ~\ref{thm1}, we find 
first order Vassiliev invariants of generic maps of $3$-manifolds
into $\R^4$, and of generic maps $S^{4k-1}\to\R^{6k-2}$, $k>1$, see
Remark ~\ref{vassinv}.  

The techniques used to prove Theorem ~\ref{thm1} allow us also to 
find restrictions on self intersections of immersions: 
\begin{thm}\label{thm3}
Let $V^{4k-1}$, be a $2$-connected closed oriented manifold
and let
$f\colon V^{4k-1}\to\R^{4k+1}$ be a generic immersion. Then
there exists an integer $d$ such that $d\cdot f$ is 
a null-cobordant immersion. Let $M^{4k}$ be a compact 
oriented manifold with $\partial M^{4k}=d\cdot V^{4k-1}$. Let
$g\colon M^{4k}\to\R^{4k+2}_+$ be a generic immersion  such 
that $\partial g =d\cdot f$. Then 
\begin{equation}\label{high0}
-\la\bar p_1^k,[M^{4k},\partial M^{4k}]\ra +
(2k+1)\sharp D_{2k+1}(g) + d\cdot L_{2k}(f) = 0.  
\end{equation}
\end{thm}
The notion $d\cdot f$ means the connected sum of $d$ copies of $f$.
The term $L_{2k}(\partial g)$ measures the linking of the 
$2k$-fold self intersection set of $\partial g$ with the rest of
its image, the term $\sharp D_{2k+1}(g)$ is the algebraic number of
$(2k+1)$-fold self intersection points of $g$, and $\bar p_1$ is the
square of the relative normal Euler class. For these notions, see
Section ~\ref{not3}. Theorem ~\ref{thm3} is proved in Section
~\ref{pfthm3}. 

For immersions $S^{4k+1}\to\R^{4k+3}$ there is a $\mod 2$-version of
the Formula ~\eqref{Shigh}. If $k$ odd then there is only one regular
homotopy class and the corresponding formula 
always vanishes, see Proposition ~\ref{prpsist}. If $k$ is even then there
are two regular homotopy 
classes and we ask whether or not in these cases the $\mod 2$-version
gives the Smale invariant, see Question ~\ref{cnj1}.  

\section{Constructions and definitions}
In this section we define all the terms appearing in the theorems
stated in the Introduction. These terms are all numerical
characteristics naturally associated to generic maps. 

\subsection{Notation and basic definitions}\label{basnot} 
We shall work in the differential category and, unless otherwise
stated, all manifolds and maps are assumed to be smooth.

It will be convenient to have a notion for the closure of a punctured
manifold: 
\begin{dfn}
If $X^n$ is an $n$-dimensional manifold with spherical boundary 
$\partial X^n\approx S^{n-1}$, then let $\hat X^n$ denote the closed
manifold obtained by gluing an $n$-disk to $X^n$ along $\partial X^n$.
\end{dfn} 

Recall that a map $f\colon M^m\to N^n$ from one manifold into another
is called {\em stable} if there exists a neighborhood $U(f)$ of $f$ in
the space $C^{\infty}(M^m,N^n)$ of smooth maps from $M^m$ to $N^n$
with the following property: For any $g\in U(f)$ there exists
diffeomorphisms $h$ of the source and $k$ of the target such that
$g=k\circ g\circ h$. 

This paper concerns the so called ``nice dimensions'' of
Mather, see ~\cite{Ma}, where the set of stable maps is open and
dense in the space of all maps. 

The notion {\em generic map} will be used throughout the paper. In general,
generic maps constitute an open dense subset of the space of all
maps. In this paper, {\em generic shall mean stable}. The requirement
that a map is generic imposes certain conditions on the following sets
associated with the map:

\begin{dfn}
If $g\colon X^n\to\R^{n+k}$ is a map of a manifold, then
\begin{itemize}
\item[{\rm (a)}]the subset $\tilde\Sigma(g)\subset X^n$ is defined as
$$
\tilde\Sigma(g)=\left\{p\in X^n\colon \rank(dg_p)\le n-1
\right\},
$$
and $\Sigma(g)=g(\tilde\Sigma(g))\subset\R^{n+k}$.
\item[{\rm (b)}] the subset $D(g)\subset\R^{n+k}$ is defined as
$$
D(g)=\left\{q\in \R^{n+k}\colon |g^{-1}(q)|\ge 2
\right\},
$$
where $|A|$ denotes the cardinality of the set $A$, and 
$\tilde D(g)=g^{-1}(D(g))\subset X^n$.
\item[{\rm (c)}] for each integer $i\ge 2$ the subset
$D_i(g)\subset\R^{n+k}$ is 
defined as 
$$
D_i(g)=\left\{q\in \R^{n+k}\colon |g^{-1}(q)|=i
\right\},
$$
and $\tilde D_i(g)=g^{-1}(D_i(g))\subset X^n$.
\end{itemize} 
\end{dfn}

In our formulas there appear certain Pontryagin numbers. 
We establish notation for these:
\begin{dfn}
If $X^{4k}$ is a closed oriented $4k$-dimensional manifold then let
$\bar p_k[X^{4k}]$ denote the Pontryagin number of $X^{4k}$ which is
associated to its $k^{\rm th}$ normal Pontryagin class.
\end{dfn}

\subsection{A triple point invariant of generic immersions}\label{nottrip}
Let $V^{4k-1}$ be a closed oriented manifold such that  
$H_{2k}(V^{4k-1};\Z)=0=H_{2k-1}(V^{4k-1};\Z)$. Let\linebreak 
$f\colon V^{4k-1}\to\R^{6k-1}$ be a generic (self-transverse)
immersion. Then $D(f)=D_2(f)\subset\R^{6k-1}$ is an embedded
$(2k-1)$-dimensional submanifold and there is an induced orientation on
$D(f)$ (since the codimension is even). 

The normal bundle of $f$ has dimension $2k$ and therefore it
admits a nonzero section $v$ over $\tilde D(f)$. 
Let $E_0$ denote the total space of the bundle of nonzero vectors in
the normal bundle of $f$.
The homology assumptions on $V^{4k-1}$ imply that
$H_{2k-1}(E_0;\Z)=\Z$. Let $[v]\in\Z$ be the homology class of 
$v(\tilde D(f))$ in $E_0$.

For $p\in D(f)$, define $w(p)=v(p_1)+v(p_2)$, where $f(p_1)=f(p_2)=p$.
Then $w$ is a normal vector field of $D(f)$ in $\R^{6k-1}$. Let $D'_v(f)$
be a copy of $D(f)$ shifted slightly along $w$. Then 
$D'_v(f)\cap f(V^{4k-1})=\emptyset$. 

\begin{dfn}\label{dfnL}
{\rm (See ~\cite{E}, ~\cite{E3})}
Let $f\colon V^{4k-1}\to\R^{6k-1}$ be an immersion as above. Define 
$$
L(f)=\lk(D'_v(f),f(S^{4k-1}))-[v],
$$
where the linking number $\lk$ is computed in $\R^{6k-1}$.
\end{dfn}
The integer $L$ as defined in Definition ~\ref{dfnL} is well
defined. That is, $L$ is independent of the choice of $v$.

\subsection{Generic maps $M^{4k}\to\R^{6k}$, linking numbers, and
triple points} \label{not2}
\begin{rmk}\label{gen2}
Let $M^{4k}$ be a compact manifold of dimension
$4k$. If $g\colon M^{4k}\to \R^{6k}$ is a {\em generic} map 
then it has the following properties:
\begin{itemize}
\item[{\rm (a)}] $D(g)=D_2(g)\cup D_3(g)$.
\item[{\rm (b)}] $\tilde D(g)\cap \tilde\Sigma(g)=\emptyset$.
\item[{\rm (c)}] At a point in $D(g)$ the self intersection is in
general position. 
\end{itemize}
\end{rmk}

If $p$ is a triple point of a generic map $g\colon M^{4k}\to\R^{6k}$
of an {\em oriented} manifold ($p\in D_3(g)$), then Remark
~\ref{gen2} (c) says that the three sheets of $M^{4k}$ meeting at $p$
intersect in general position. Therefore, the tangent space
$T_p\R^{6k}$ splits into a direct sum of the three oriented
$2k$-dimensional normal spaces of the  
sheets. Thus, there is an orientation induced on each triple point of
$g$. 
\begin{dfn}\label{dfnt}
Let $g\colon M^{4k}\to\R^{6k}$ be a generic smooth map. Define 
$t(g)$ as the algebraic number of triple points of the map $g$.
\end{dfn} 

For the sake of the next definition we separate the cases into:
\begin{itemize}
\item[(a)] The manifold $M^{4k}$ is closed.
\item[(b)] The manifold $M^{4k}$ has non-empty boundary. 
\end{itemize} 
In case (a), let $g\colon M^{4k}\to \R^{6k}$ be any
generic map. 

Let $\R^{6k}_+$ denote a closed half-space of $\R^{6k}$.
In case (b), let $g\colon M^{4k}\to\R^{6k}_+$ 
be any generic map such that its restriction $\partial g$ to the
boundary  $\partial M^{4k}$ is an immersion and such that
$g^{-1}(\partial\R^{6k}_+)=\partial M^{4k}$.  

Then $D(g)$ is an immersed $2k$-dimensional submanifold of $\R^{6k}$
with non-generic triple self intersections at the triple points of $g$
and $\Sigma(g)$ is an embedded $(2k-1)$-dimensional 
manifold. In case (a) $\Sigma(g)$ is the boundary of $D(g)$,  
and in case (b) $\Sigma(g)$ is a part of the
boundary of $D(g)$ (the other part is the double points of $\partial g$).

Since $\dim(\R^{6k})-\dim(M^{4k})=2k$ is an even number, there is an
induced orientation on $D(g)$, which in turn induces an orientation of
$\Sigma(g)$.   

Let $\Sigma'(g)$ be a copy of $\Sigma(g)$ shifted slightly along the outward
normal vector field of $\Sigma(g)$ in $D(g)$. Then 
$\Sigma'(g)\cap g(M^{4k})=\emptyset$.

\begin{dfn}\label{dfnl}
{\rm (See ~\cite{S})} 
Let $g$ be a map as above. Define $l(g)$ as the linking number of
$g(M^{4k})$ and $\Sigma'(g)$, in $\R^{6k}$ in case {\rm (a)}, and  
in $(\R^{6k}_+,\partial\R^{6k}_+)$ in case {\rm (b)}.   
\end{dfn}

\subsection{Generic maps $M^{4k}\to\R^{6k-1}$, cusps, and Euler
classes}\label{not1} 
\begin{rmk}\label{gen1}
Let $M^{4k}$ be a compact manifold of dimension
$4k$. If $g\colon M^{4k}\to \R^{6k-1}$ is a {\em generic} map then it has
the following properties:
\begin{itemize}
\item[{\rm (a)}] For $k>1$, $D(g)=D_2(g)\cup D_3(g)$.
For $k=1$, $D(g)=D_2(g)\cup\dots \cup D_5(g)$.  
\item[{\rm (b)}] For $k>1$, 
$\tilde D_3(g)\cap\tilde\Sigma(g)=\emptyset$. 
For $k=1$, $(\tilde D_4(g)\cup\tilde D_5(g))\cap\tilde \Sigma(g)=\emptyset$.
\item[{\rm (c)}] $\tilde\Sigma(g)$ is a $2k$-dimensional submanifold
of $M^{4k}$. At each point $p\in\tilde\Sigma(g)$, $\rank(dg)=4k-1$. Thus,
the kernel $\ker(dg)$ of $dg$ is a
$1$-dimensional subbundle of the restriction $T M^{4k}|\tilde\Sigma(g)$
of the tangent bundle to $\tilde\Sigma(g)$. 

Moreover, the fiber $\ker(dg)_p\subset T_p M^{4k}$ does not
lie in $T_p\tilde\Sigma(g)\subset T_p M^{4k}$ for all but finitely  
many $p\in\tilde\Sigma(g)$. The finitely many exceptional points are called
{\em $\Sigma^{1,1}$-points} or {\em cusps}.  
\item[{\rm (d)}] No cusp is in $\tilde D(g)$. At a point in $D(g)-\Sigma(g)$
the self intersection is in general position. At a point in
$D(g)\cap\Sigma(g)$ the smooth sheet (or sheets if $k=1$) of
$g(M^{4k})$ meets $\Sigma(g)$ in general position. 
(Since no cusp is in $\tilde D(g)$, the tangent space of
$\Sigma(g)$ is well defined at all points in $D(g)\cap\Sigma(g)$).  
\end{itemize}
\end{rmk}

It is proved in ~\cite{S2}, Appendix 1, that if $q$ is a
$\Sigma^{1,1}$-point of a generic map $g\colon M^{4k}\to\R^{6k-1}$ of
an {\em oriented} manifold then there are induced orientations on
$T_{g(q)}\R^{6k-1}$ and $T_q M^{4k}$. 
Taking the product of these orientations, 
a sign is associated to each
$\Sigma^{1,1}$-point. 
\begin{dfn}\label{dfnS11}
If $g\colon M^{4k}\to\R^{6k-1}$ is a generic map of an oriented
manifold then let $\sharp\Sigma^{1,1}(g)$ denote its algebraic number
of $\Sigma^{1,1}$-points. 
\end{dfn}

\begin{dfn}\label{dfnxi}
If $g\colon M^{4k}\to\R^{6k-1}$ is a generic map then let 
$\xi(g)$ denote the $2k$-dimensional vector bundle over
$\tilde\Sigma(g)$, the fiber of which over 
$p\in\tilde\Sigma(g)$ is 
$\xi(g)_p=T_{g(p)}\R^{6k-1}/dg(T_pM^{4k})$.
\end{dfn}
The total space $E(g)$ of $\xi(g)$ is orientable (and even
oriented), see Lemma ~\ref{orxi}. In this situation, the Euler number
of $\xi(g)$ is a well-defined integer (in the case when
$\tilde\Sigma(g)$ is non-orientable it is necessary to use homology
with twisted coefficients, see Example ~\ref{ex1} and Section
~\ref{twxi}.) 

\begin{dfn}      
If $g\colon M^{4k}\to\R^{6k-1}$ is a generic map of an oriented
manifold then let $e(\xi(g))$ denote the Euler number of the bundle
$\xi(g)$ over $\tilde\Sigma(g)$ (see Definition ~\ref{dfnxi}).
\end{dfn}

\subsection{Self intersection points of high multiplicity}\label{not3}
Let $M^{4k}$ be a compact oriented manifold of dimension $4k$ and let
$g\colon M^{4k}\to\R^{4k+2}$ be a generic (self transverse)
immersion. Then $g$ has isolated $(2k+1)$-fold self intersection
points. At such a point the 
sheets of $M^{4k}$ meet in general position. The oriented normal
spaces to the sheets of $M^{4k}$ at $p$ induce an orientation of
$T_p\R^{4k+2}$.
\begin{dfn}
For $g$ as above let $\sharp D_{2k+1}(g)$ denote the algebraic number of
$(2k+1)$-fold self intersection points of $g$.
\end{dfn}

Let $V^{4k-1}$ be a $2$-connected closed oriented manifold and let 
$f\colon V^{4k-1}\to\R^{4k+1}$ be a generic immersion. Then the
$2k$-fold self intersection points of $f$ form a closed $1$-manifold
$D_{2k}(f)\subset\R^{4k+1}$. Since $V^{4k-1}$ is $2$-connected it has
a (unique up to homotopy) normal framing $(n_1,n_2)$. Let $D'_{2k}(f)$
be the manifold which is obtained when $D_{2k}(f)$ is shifted slightly
along the vector field, $p\mapsto n_1(p_1)+\dots+n_1(p_{2k})$, for
$p\in D_{2k}(f)$, $p=f(p_1)=\dots=f(p_{2k})$.
\begin{dfn}
For an immersion $f$ as above, define
$$
L_{2k}(f)=\lk(f(S^{4k-1}),D'_{2k}(f)),
$$
where the linking number $\lk$ is computed in $\R^{4k+1}$.
\end{dfn}
Let $f\colon V^{4k-1}\to\R^{4k+1}$ be an immersion. Then there is a
nonzero integer $d$ such that $d\cdot f=f\sharp\dots\sharp f$ 
bounds an immersion $g\colon M^{4k}\to\R^{4k+2}_+$. (See Section
~\ref{mpscob} for this fact and ~\cite{E3} for the connected sum
operation on {\em generic} immersions.)

The normal bundle $\nu_g$ of $g$ is trivial over 
$d\cdot V^{4k-1}=\partial M^{4k}$, since $V^{4k-1}$ is
$2$-connected. Moreover, its trivialization is homotopically
unique. Hence, the Euler class $\bar e$ of $\nu_g$ can be considered
as a class in $H^{2}(M^{4k},\partial M^{4k};\Z)$. We introduce the
following notation: 
$$
\bar p_1=p_1(\nu_g)=\bar e^2.
$$  
Finally, we let $\langle\bar p_1^k,[M^{4k},\partial M^{4k}]\rangle$ denote the
evaluation of $\bar p_1^k$ on the orientation class 
$[M^{4k},\partial M^{4k}]$.

\section{Some remarks on the main results}
In this section we make several remarks concerning
extensions of our formulas and concerning their relations
to other results.

\subsection{The Smale invariant formulas}  
\begin{rmk}\label{3to5} 
Consider the first case, $k=1$, in Theorem ~\ref{thm1}. 
We can rewrite Formulas ~\eqref{grom1}, ~\eqref{grom2}, and  ~\eqref{Shigh}
using the following well-known identities for an oriented closed
$4$-manifold $X$,  
$$
-\bar p_1[X]=p_1[X]=3\sigma(X), 
$$
where $p_1[X]$ is the Pontryagin class of the tangent bundle evaluated
on the orientation class of $X$ and $\sigma(X)$ denotes the signature
of $X$. 

The signature of a non-closed $4$-manifold $M^4$ can still be
defined, see ~\cite{Ki}. In the general case the quadratic form on 
$H_2(M^4;\Z)$ may be degenerate but if $M^4$ has spherical boundary
then $\sigma(M^4)=\sigma(\hat M^4)$. 

The resulting formulas are   
\begin{align}\label{grom351}
\Omega (f) &= \frac12\left(3\sigma(M^4) +\sharp\Sigma^{1,1}(g)\right),\\ 
&= \frac12\left(3\sigma(M^4) +e(\xi(g))\right)\label{grom352}
\end{align}
replacing Equations ~\eqref{grom1} and ~\eqref{grom2}, respectively, and
\begin{equation}\label{S35}
\Omega (f) = \frac12\left(3\sigma(M^4) + 3 t(g) - 3 l(g) +L(f)\right),
\end{equation}
replacing Equation ~\eqref{Shigh}.
\end{rmk}
\begin{rmk}\label{vassinv}
Formulas ~\eqref{grom351} and ~\eqref{grom352} give invariants of
arbitrary generic maps of closed oriented $3$-manifolds $V^3$ into
$\R^4$: 

Fix some oriented manifold $M^4$ such that $\partial M^4=V^3$.
Let $f\colon V^3\to\R^4$ be any generic map. Then there exists a
generic map $g\colon M^{4}\to\R^6_+$ such that
$\partial g=f$. Define 
$$
v(f)=\sigma(M^4)+\sharp\Sigma^{1,1}(g)=\sigma(M^4)+ e(\xi(g)).
$$ 
By Lemmas ~\ref{s11clo} and ~\ref{es11}, $v$ is well defined. Moreover,
$v$ does not change during homotopies through generic maps. More
precisely,  $v$ does not change under homotopies
avoiding maps with cusps and changes by $\pm 1$ at cusp
instances. Thus, $v$ is a Vassiliev invariant of degree one.

In the same way we can define first order Vassiliev invariants of
generic maps $S^{4k-1}\to\R^{6k-2}$: Let $f$ be such a map. Pick a
generic map $F\colon D^{4k}\to\R^{6k-1}_+$ extending $f$. Define
$v(f)=\sharp\Sigma^{1,1}(F)$. Then $v$ is well defined by Lemma
~\ref{s11clo}. Again $v$ is an invariant of generic
maps and changes by $\pm 1$ at cusp instances. 
\end{rmk}

\begin{rmk}\label{ext35}
The expressions for $\Omega (f)$ given in Formulas 
~\eqref{grom351}, ~\eqref{grom352}, and ~\eqref{S35} make sense also when 
$f\colon V^3 \to\R^5$ is a framed generic immersion of any oriented
$3$-manifold $V^3$: 

The signature still makes sense, see Remark ~\ref{3to5}. The
definition of $L(f)$ can be extended to this case as follows:   

Let $n_1,n_2$ be the two linearly independent normal vectors in the
normal bundle of $f$, which give its framing. The fiber of the normal
bundle of $f$ at a point $p\in V^3$ can be identified with
$T_{f(p)}\R^5/df(T_pV^3)$. Using this identification we let 
$w$ be the vector field along $D(f)$ defined by
$w(q)=n_1(p_1)+n_1(p_2)$, where $f(p_1)=f(p_2)=q$.
Let $D'(f)$ be the manifold which results from pushing $D(f)$ a small
distance along $w$. Then $D'(f)\cap f(V^3)=\emptyset$. 

Define $L(f)$ as the linking number of
$D'(f)$ and $f(V^3)$ in $\R^5$. (In the special case of a homology
sphere this definition agrees with Definition ~\ref{dfnL}.)   

Let $\pi^{\rm s}(3)$ denote the stable homotopy group
$\pi_{N+3}(S^N)$, $N$ large. If $f\colon V^3\to\R^5$ is framed generic
immersion of any $3$-manifold, then $\Omega(f)$ (as given in
Formulas ~\eqref{grom351}, ~\eqref{grom352}, and ~\eqref{S35})
reduced modulo $24$ gives the element in 
$\pi^{\rm s}(3)=\Z_{24}$ realized by the framed immersion $f$.
\end{rmk}

\begin{rmk} 
In higher dimensions, we can find analogs of Formulas
 ~\eqref{grom351}, ~\eqref{grom352}, and ~\eqref{S35} 
when $\hat M^{4k}$ is almost parallelizable. If this is the case
then $\bar p_k[\hat M^{4k}]$ can be replaced by
the appropriate multiple of the signature:
$$ 
-\bar p_k[\hat M^{4k}] = p_k[\hat M^{4k}] = 
\frac {(2k)!}{2^{2k}(2^{2k-1} -1) B_k} \sigma(M^{4k})
$$
where $B_k$ is the $k$-th Bernoulli number, see ~\cite{MK}.
\end{rmk}
\begin{rmk} 
Formulas ~\eqref{grom1}, ~\eqref{grom2}, and  ~\eqref{Shigh} reduced 
modulo the order of the image of the $J$-homomorphism give the
element represented by the immersion $f\colon S^{4k-1}\to\R^{4k+1}$
with its homotopically unique normal framing in $\pi^{\rm s}(4k-1)$.
\end{rmk}

\begin{rmk} 
The right-hand sides of 
Equations ~\eqref{grom1}, ~\eqref{grom2}, and  ~\eqref{Shigh} can be
altered as follows:  

Replace the normal Pontryagin class of $\hat M^{4k}$
by the corresponding class of the normal bundle of $M^{4k}$ considered
as an element in the relative cohomology group 
$H^{4k}(M^{4k},\partial M^{4k};\Z)$. 
(This is possible since the normal bundle over the boundary is
trivialized; it inherits a homotopically unique normal framing from
the codimension two immersion $f$.) These altered formulas vanish for
any generic $g$ satisfying the conditions in Theorem ~\ref{thm1}.
(For $k=1$, Theorem ~\ref{thm3} is one of these altered 
formulas.) 

The proof of these facts are the same as the proofs of the formulas
themselves. 
\end{rmk}

\begin{rmk}\label{exfel}
Formula ~\eqref{grom2} was inspired by Exercise (c) on page 65 in
Gromov's book ~\cite{G}. The exercise reads:

``Let $V$ be a closed oriented $4$-manifold and let 
$f\colon V\to\R^5$ be a generic $C^\infty$-map. 
Then the singularity $\Sigma=\Sigma^1_f$ is a smooth closed surface in
$V$ such that $\rank_v df=3$ for all $v\in\Sigma$. Let
$g\colon\Sigma\to Gr_3\R^5$ be the map which assigns the image
$D_f(T_v)$ (which is a $3$-dimensional subspace in $\R^5$) to each
point $v\in V$. Prove, for properly normalized Euler
form $\omega$ in $Gr_3\R^5$, [which is a closed $SO(5)$ invariant
$2$-form on the Grassmann manifold $Gr_3\R^5=Gr_2\R^5$], the equality
$\int_\Sigma g^\ast(\omega)=p_1(V^4)$ for a natural orientation in
$\Sigma(g)$ and for the first Pontryagin number $p_1(V)$ of $V$.''

Using the notation introduced in Section ~\ref{not1}, the Euler class
$g^{\ast}(\omega)$ is the Euler class of the bundle 
$\xi(f)$. Since  $\Sigma$ ($\tilde\Sigma(f)$ in our notation) is
sometimes non-orientable, the statement
{\em ``for a natural orientation in $\Sigma$''} should be understood
in terms of homology with local coefficients, see Section ~\ref{twxi}. 

We give an example with (necessarily) non-orientable singularity surface:   
\begin{exm}\label{ex1}
Let $f\colon \C P^2\to\R^5$ be any generic map. It follows from
Lemma ~\ref{orxi} below that the determinant bundle
$\det(T\tilde\Sigma(f))$ 
of $\tilde\Sigma(f)$ is isomorphic to  $1$-dimensional vector
bundle $\ker(df)$ over $\tilde\Sigma(f)$. Thus the orientability of
$\tilde\Sigma(f)$ implies that $\ker(df)$ is a trivial line bundle. 
Hence, by Lemma (A) on page 49 in ~\cite{G}, there exists a function
$\phi$ such that $d\phi(\ker(df))\ne 0$. 
Then $F\colon\C P^2\to\R^5\times\R$,  
$F(x)=(f(x),\phi(x))$ is an immersion. However, since $p_1[\C P^2]=3$ is
not a square, $\C P^2$ does not immerse in $\R^6$. (If it did immerse,
$p_1(\C P^2)$ would equal $\bar e^2$, where $\bar e$ is the Euler
class of the normal bundle of the immersion.) This contradiction shows
that $\tilde\Sigma(f)$ cannot be orientable.
\end{exm}
\end{rmk}

\subsection{Regular homotopy of embeddings in high codimension}
\begin{rmk}
In ~\cite{HM}, Hughes and Melvin pose the following question: 
{\em ``For a given $n$, what is the largest possible value of $k$ such
that $\Emb(S^n,\R^k)\ne 0$?''}
Here, $\Emb(S^n,\R^k)$ denotes the set (group) of regular homotopy
classes which can be represented by embeddings. 
Combining Corollary ~\ref{corthm1} and the theorem of Kervaire stated
in the Introduction gives the answer for $n=4k-1$:

There are infinitely many regular homotopy classes of immersions
$S^{4k-1}\to\R^{q}$ containing embeddings if $4k+1\le q\le 6k-1$. If 
$q\ge 6k$ then there is only one such class.
\end{rmk}

\subsection{Bounding immersions}
\begin{rmk}\label{nongen}
Let $V^{4k-1}\approx S^{4k-1}$ in Theorem  ~\ref{thm3}.
Replacing the term $\langle\bar p_1^k,[M^{4k},\partial M^{4k}]\rangle$
in Formula ~\eqref{high0} by 
$\bar p_1^k[\hat M^{4k}]$, the equation would still hold for $k>1$.
However, for $k=1$, the left-hand side would (using Theorem ~\ref{thm1}
(b)) equal $d\cdot\Omega(f)$.

This difference between $k=1$ and $k>1$ comes from the 
fact that if $M^{4k}$, $k>1$ is a manifold with spherical boundary 
$\partial M^{4k}= S^{4k-1}$ and $g\colon M^{4k}\to\R^{4k+2}_+$ is an
immersion, 
then $\bar p_1(\hat M^{4k})$ corresponds to $p_1(\nu_g)$, where
$\nu_g$ is the normal bundle of $g$. In this case we can consider
$\bar p_1$ as an absolute class already in $M^{4k}$ since 
$$
H^4(\hat M^{4k};\Z)\approx H^4(M^{4k},\partial M^{4k};\Z)\approx
H^4(M^{4k};\Z). 
$$
This is not the case if $k=1$, then $H^4(M^4;\Z)=0$.
\end{rmk}
\begin{rmk}
According to Remark ~\ref{nongen}, Theorem ~\ref{thm3} should not be
considered as a generalization of Theorem ~\ref{thm1} (b). Instead it can
be considered as a generalization of a theorem in ~\cite{CS}, in which
the number of triple points of a surface immersed in the upper
half-space and bounding a given regular plane curve has been computed.
\end{rmk}

\section{Smale invariants of embeddings $S^{4k-1}\to\R^{4k+1}$}
In this section we shall state a result which will be used in 
our proof of Theorem \ref{thm1}. The result is due to Hughes and
Melvin ~\cite{HM}. 

Let $f\colon S^{4k-1}\to\R^{4k+1}$ be a smooth embedding. It is a
well known fact that any such $f$ admits a Seifert-surface. That is,
there exists a smooth orientable $4k$-manifold $M^{4k}$, with
spherical boundary and an embedding $g\colon M^{4k}\to\R^{4k+1}$ such
that the restriction $\partial g$ of $g$ to $\partial M^{4k}$ equals
$f$.

Since $M^{4k}$ immerses into $\R^{4k+1}$ it follows that $\hat M^{4k}$
is almost parallelizable. The immersion $g$ gives a framing of the
stable tangent bundle of $\hat M^{4k}$ along $M^{4k}$. The obstruction
to extending this trivialization over the added disk can on the one
hand be expressed in terms of the $k^{\rm th}$ Pontryagin class of
$\hat M^{4k}$ (see ~\cite{MK}) and on the other hand it can be
expressed in terms of 
the Smale invariant of $f$. Using this observation one can prove the
following:
\begin{lma}\label{regemb}
There are embeddings $S^{4k-1}\to\R^{4k+1}$ which are not regularly
homotopic to the standard embedding. Moreover, if $f\colon
S^{4k-1}\to\R^{4k+1}$ is an embedding and 
$g\colon M^{4k}\to\R^{4k+1}$ is a Seifert-surface of $f$ then 
$$
\Omega(f)=\frac1{a_k(2k-1)!}p_k[\hat M^{4k}],
$$
where $\Omega(f)$ is the Smale invariant of $f$, $a_k=2$ for $k$ odd,
and $a_k=1$ for $k$ even.
\end{lma}

\begin{pf}
As mentioned above, this is proved in ~\cite{HM}.
Here we present a simple, slightly different proof. 
Let $\Imm(S^{4k-1},\R^{4k+1})$ denote the set of regular homotopy
classes of  
immersions $S^{4k-1}\to\R^{4k+1}$ and let $\pi_{4k-1}(SO)\cong\Z$
denote the stable 
homotopy group of the orthogonal group. 

There is a map $\Imm(S^{4k-1},\R^{4k+1})\to\pi_{4k-1}(SO)$: 
Given an immersion, lift it to 
$\R^N$, $N$ large. This immersion has a homotopically unique (normal)
framing. (The unique framing of the immersion in $\R^{4k+1}$ plus the
trivial framing of $\R^{4k+1}$ in $\R^N$.) Deform this lifted
immersion to the standard embedding and compare the 
induced framing with the standard one (obtained from the standard
embedding of $S^{4k-1}\subset\R^{4k}\subset\R^{4k+1}\subset\R^N$).

This map has an inverse: Via the Hirsch lemma, any given framing
on the standard embedding $S^{4k-1}\to\R^N$ can be used to push the framed
immersion down to an immersion $S^{4k-1}\to\R^{4k+1}$, the regular homotopy
class of which is well-defined.

It is well known that $\pi_{4k-1}(SO)=\Vect(S^{4k})$, where
$\Vect(S^{4k})$ denotes the group of stable equivalence classes of
vector bundles on $S^{4k}$. If $\eta\in\Vect(S^{4k})$ is a stable
bundle then let $[\eta]$ denote the corresponding element in
$\pi_{4k-1}(SO)$ and if $f\colon S^{4k-1}\to\R^{4k+1}$ is an immersion
then let $\eta_f$ denote the stable bundle corresponding to $f$. 

Let $p_k\colon\Vect(S^{4k})\to\Z$ denote the map defined by
$\eta\mapsto\langle p_k(\eta),[S^{4k}]\rangle$. 
Lemma 2 in ~\cite{MK} says that 
$\eta\mapsto\langle p_k(\eta),[S^{4k}]\rangle=a_k(2k-1)![\eta]$, where
$a_k=2$ if $k$ is odd and $a_k=1$ if $k$ is even.

Consider the map $\Imm(S^{4k-1},\R^{4k+1})\to\Z$ defined as the
composition of 
the above two, i.e. $f\mapsto p_k(\eta_f)$. Since the first map is onto,
this composition is a map of $\Z$ onto $a_k(2k-1)!\cdot\Z$. 
Hence, it is $\pm a_k(2k-1)!\cdot\Omega$, where $\Omega$ is the Smale
invariant. That is, 
$\pm a_k(2k-1)!\cdot\Omega(f)=\langle p_k(\eta_f),[S^{4k}]\rangle$.

Consider the standard embedding of $S^{4k-1}$ in $\R^N$ with framing
induced from an immersion $f\colon S^{4k-1}\to\R^{4k+1}$. Assume that it
bounds a framed immersion $F\colon M^{4k}\to\R^{N+1}_+$. Let $\nu$ denote
the stable normal bundle of the closed manifold $\hat M^{4k}$. 
It is straightforward to find a map $g\colon\hat M^{4k}\to S^{4k}$ of
degree $1$ such that $g^\ast(\eta_f)=\nu$. Thus,
$$
\langle p_k(\nu),[\hat M^{4k}]\rangle=
\langle p_k(\eta_f),[S^{4k}]\rangle.
$$ 

Note that $p_k[\hat M^4]=-p_k(\nu)$ (since $\hat M^4$ is almost
parallelizable). 
Hence,\linebreak
$p_k[\hat M^{4k}]=\pm a_k(2k-1)!\cdot\Omega(f)$. 

Especially, one can apply this when $f$ is an embedding and $M^{4k}$ is a
Seifert surface of $f$. This proves the formula in Lemma
~\ref{regemb}. 

The existence statement is proved by noting that there exist almost
parallelizable $4k$-manifolds $\hat X^{4k}$, with $p_k[\hat X^{4k}]\ne
0$, and with a handle decomposition consisting of one $0$-handle,
$2k$-handles, and one $4k$-handle. 
Let $X^{4k}$ be a punctured $\hat X^{4k}$. Then $X^{4k}$ immerses into
$\R^{4k+1}$ by the Hirsch lemma. Its $2k$-skeleton actually embeds in
$\R^{4k+1}$ by 
general position, but $X^{4k}$ is a regular neighborhood of this
$2k$-skeleton in $\hat X^{4k}$ so $X^{4k}$ also embeds. The boundary
of an embedded $X^{4k}$ is an embedded sphere with non-trivial Smale
invariant.      
\end{pf}

\section{Generic $4k$-manifolds in $6k$-space and \\ 
the proof of Theorem ~\ref{thm1} {\rm (b)}}
In this section we present some lemmas on generic maps and use
these to prove Theorem ~\ref{thm1} (b).

\subsection{Generic $4k$-manifolds in $6k$-space}
We shall show that the difference of triple points and singularity linking
of a generic map of a closed oriented $4k$-manifold in $6k$-space is a
cobordism invariant. The proof will use the following well known fact
about linking numbers.
\begin{lma}\label{lkint}
Let $X$ and $Y$ be relative oriented cycles of complementary
dimensions in $\R^n\times I$, where $I$ is the unit interval. Let
$\partial_i X$ and $\partial_i Y$ be 
$X\cap\R^{n}\times\{i\}$ and $Y\cap\R^{n}\times\{i\}$, respectively,
for $i=0,1$. Assume that $\partial_i X\cap\partial_i Y=\emptyset$ for
$i=0,1$ and that the intersection of $X$ and $Y$ is transverse. Then
$$
X\bullet Y=\lk(\partial_1 X,\partial_1 Y)-\lk(\partial_0 X,\partial_0 Y),
$$ 
where $X\bullet Y$ denotes  the intersection number of $X$ and $Y$.
\qed
\end{lma}
\begin{lma}\label{cobinv}
Let $M^{4k}_0$ and $M^{4k}_1$ be two closed and oriented manifolds. Let 
$k_i\colon M^{4k}_i\to\R^{6k}$, $i=0,1$, be generic smooth maps. Assume
that there exists an oriented cobordism 
$K\colon W^{4k+1}\to\R^{6k}\times I$, 
where $I$ is the unit interval, joining $k_0$ to $k_1$. Then
$$
t(k_0)-l(k_0)=t(k_1)-l(k_1).
$$ 
(For notation, see Definitions ~\ref{dfnt} and ~\ref{dfnl}.)
\end{lma}
\begin{pf}
We may assume that the map $K$ is generic. Then $K$ does not have any
$4$-fold self intersection points. The triple points of $K$ form a
$1$-manifold $D_3(K)$ and there is an induced orientation of
$D_3(K)$. The boundary of $D_3(K)$ consists of three types of points:.
\begin{itemize}
\item[(a)] Triple points of $k_1$,
\item[(b)] triple points of $k_0$, and
\item[(c)] those double points of $K$ which are also singular values. 
(That is, points in $\Sigma(K)\cap D(K)$.)  
At such a double point a nonsingular sheet of $W^{4k+1}$ meets the
$2k$-dimensional singular locus $\Sigma(K)$ of $K$ transversely. (In other
words, it 
is a double point of type $\Sigma^{1,0}+\Sigma^0$.)  
\end{itemize}
There is an induced orientation on the manifold
$D_2(K)$ of nonsingular double points of $K$. Therefore,
there is an induced orientation on its boundary $\Sigma(K)$. Hence, there
are also induced orientations on the points of type (c) above. 

Note that a point $p$ of type (a) is a positive triple point of $k_1$ if
and only if the orientation of $D_3(K)$ close to $p$ points towards
$p$ (out of $D_3(K)$), that a point $q$ of type (b) is a positive
triple point of $k_0$ if and only if the orientation points away from
$q$ (into $D_3(K)$), and that a point $r$ of type (c) has positive
sign if and only if the orientation points towards $r$.

Let $A(K)$ denote the algebraic number of points of type (c). Since
the triple curves of $K$ give a cobordism between the points of type
(c) and those of type (a) or (b), it follows that
\begin{equation}\label{trip}
A(K)=t(k_1)-t(k_0).
\end{equation} 

Now consider the $2k$-dimensional manifold $\Sigma(K)$ of singular values
of $K$. The boundary of $\Sigma(K)$ consists of $\Sigma(k_1)$ and
$\Sigma(k_2)$. Let $n$ be the outward normal vector
field of $\Sigma(K)$ in $D_2(K)$. Pushing $\Sigma(K)$ along $n$ we
obtain an oriented  manifold $\Sigma'(K)$. 

Close to each  
point $p$ of type (c) above, $\Sigma'(K)$ intersects $K(W^{4k+1})$ at one
point with local intersection number equal to the sign of $p$,
$\Sigma'(K)$ does not intersect $K(W^{4k+1})$ in other points and the
boundary of 
$\Sigma'(K)$ is $\Sigma'(k_1)-\Sigma'(k_0)$. Therefore, by Lemma
~\ref{lkint}, 
$$
A(K)=l(k_1)-l(k_0).
$$
\end{pf}

\begin{lma}\label{lma1}
Let $M^{4k}$ be a closed and oriented manifold and let 
$h\colon M^{4k}\to\R^{6k}$ be a generic smooth map. Then
$$
-\bar p_k[M^{4k}]+3t(h)-3l(h)=0.
$$ 
(For notation, see Definitions ~\ref{dfnt} and ~\ref{dfnl}.)
\end{lma}
\begin{pf}
For immersions $h$, $l(h)=0$ and Lemma ~\ref{lma1} is a  theorem of
Herbert, see ~\cite{He}. 

Let $\Imm^{SO}(4k,2k)$ denote the cobordism group of immersions of oriented
$4k$-dimensional manifolds in $\R^{6k}$. Let $\Omega_{4k}(\R^{6k})$ be
the cobordism group of generic maps of oriented $4k$-manifolds into
$\R^{6k}$. Then, of course, $\Omega_{4k}(\R^{6k})\approx\Omega_{4k}$,
where $\Omega_{4k}$ is the cobordism group of oriented
$4k$-manifolds. 

A theorem of Burlet ~\cite{Bu} says that the cokernel of the
natural map $\Imm^{SO}(4k,2k)\to\Omega_{4k}(\R^{6k})$ is
{\em finite}. (Sketch of proof of Burlet's theorem:
Note that 
$\Imm^{SO}(4k,2k)\approx
\pi^{\rm s}_{6k}(MSO(2k))\approx\pi_{6k+K}(\Sigma^KMSO(2k))$ 
and that 
$\Omega_{4k}\approx\pi_{4k+K}(MSO(K))\approx\pi_{6k+K}(MSO(2k+K))$ for $K$
sufficiently large. 
Apply
Serre's theorem, saying that the rational stable Hurewicz
homomorphism $\pi^{\rm s}_i(X)\otimes\Q\to H_i(X;\Q)$ is an isomorphism for
any space $X$, pass to (co)homology, and use the Thom isomorphism and
the well known ring $H^*(BSO(m);\Q)$.)  

Now, by Lemma ~\ref{cobinv}, the difference $3t(g)-3l(g)$ is invariant
under cobordism  
and hence gives rise to a homomorphism $\Omega_{4k}(\R^{6k})\to\Z$. Denote it
$\Lambda$. By Herbert's theorem $\Lambda$ equals $\bar p_k$ on the
image of the group $\Imm^{SO}(4k,2k)$ in $\Omega_{4k}(\R^{6k})$. Since this
subgroup has finite index, $\Lambda$ and $\bar p_k$ agree on the whole group.
\end{pf}
\subsection{Proof of Theorem ~\ref{thm1} (b)}\label{pfthm1b}
We show first that if $f\colon S^{4k-1}\to\R^{4k+1}$ is a generic
immersion then the expression
\begin{equation}\label{wdexpr}
\Theta(f)=-\bar p_k[\hat M^{4k}]+3t(g)-3l(g)+L(\partial g),
\end{equation}
is independent of the choice of the map $g\colon M^{4k}\to\R^{6k}_+$ and
invariant under regular homotopy of $f$.

Let $g_0\colon M^{4k}_0\to\R^{6k}_+$ and 
$g_1\colon M^{4k}_1\to\R^{6k}_+$ be two generic maps with 
$\partial g_0$ and $\partial g_1$ both regularly homotopic to $j\circ
f$.  

Then $\partial g_0$ is regularly homotopic to $\partial g_1$. Let 
$k_t$ be a generic regular homotopy between them and let
$K\colon S^{4k-1}\times I\to\R^{6k-1}\times I$ be the map 
$K(x,t)=(k_t(x),t)$. To each triple point instance of $k_t$ there
corresponds a change in $L(k_t)$ by $\pm 3$ (see ~\cite{E} or
~\cite{E3}), and also a triple point of $K$ which has the same sign as
the sign of the 
change in $L(k_t)$. Thus,
\begin{equation}\label{dL}
L(\partial g_1)-L(\partial g_0)=3t(K).
\end{equation}
Let $-g_0\colon -M^{4k}_0\to\R^{6k}_-$ be the map $r\circ g_0$, where
$r\colon\R^{6k}\to\R^{6k}$ is the reflection in $\partial \R^{6k}_+$. 
Using $K$ to
glue $-g_0$ and a vertically translated copy of $g_1$, we get a map
$h\colon -\hat M^{4k}_0\sharp \hat M^{4k}_1\to\R^{6k}$. Then, by Lemma
~\ref{lma1},  
\begin{align*}
0&=-\bar p_k[-\hat M^{4k}_0\sharp \hat M^{4k}_1]+3t(h)-3l(h)\\
&=-\bar p_k[\hat M^{4k}_1]+\bar p_k[\hat M^{4k}_0]
+3\left(t(g_1)-t(g_0)+t(K)\right)-3\left(l(g_1)-l(g_0)\right).
\end{align*}
Together with Equation ~\eqref{dL}, this implies that $\Theta(f)$ is
independent of the choice of $g$ and also that $\Theta$ is invariant
under regular homotopy. (Actually, it implies something stronger:
$\Theta$ only depends on the regular homotopy class of $j\circ f$ in
$\R^{6k-1}$. This observation leads to Corollary ~\ref{corthm1}.) 

Clearly, $\Theta$ is additive under connected sum and hence induces a
homomorphism $\Theta\colon\Imm(S^{4k-1},\R^{4k+1})\to\Z$. Thus,
$\Theta=c\cdot\Omega$, where $c\in\Z$ is some constant and $\Omega$ is the
Smale invariant. 

Then according to Lemma ~\ref{regemb},  $c={a_k(2k-1)!}$.
Indeed, let $g\colon M^{4k}\to\R^{4k+1}\subset\R^{6k-1}$ be a
Seifert-surface of 
an embedding $\partial g\colon S^{4k-1}\to\R^{4k+1}$. Push the
interior of $g(M^{4k})$ into $\R^{6k}_+$. Then  
$l(g)=0$ and $L(\partial g)=0$ since $g$ has neither singularities
nor self intersections. Hence, $\Theta(f)=-\bar p_k[\hat M^{4k}]$ and the
latter is ${a_k(2k-1)!}\cdot\Omega(f)$ by Lemma ~\ref{regemb}.
\qed

\section{Generic $4k$-manifolds in $(6k-1)$-space and\\ 
the proof of Theorem ~\ref{thm1} {\rm (a)}}
In this section we show that the right-hand sides in
~\eqref{grom1} and ~\eqref{grom2} agree and demonstrate, after
introducing a small amount of cobordism theory, how 
Theorem ~\ref{thm1} (a) follows from Theorem ~\ref{thm1} (b).
\subsection{A vector bundle over singularities}\label{twxi}
Let $g\colon M^{4k}\to\R^{6k-1}$ be a generic map of a compact
orientable manifold (see Remark ~\ref{gen1}). Then
the singularity set $\tilde\Sigma(g)$ of the map $g\colon
M^{4k}\to\R^{6k-1}$  is a $2k$-dimensional submanifold  
of $M^{4k}$, which is not necessarily orientable, see Example
~\ref{ex1}. Recall, see Definition ~\ref{dfnxi}, 
that $\xi(g)$ is the $2k$-dimensional vector bundle over
$\tilde\Sigma(g)$, the fiber of which at a point $p\in\tilde\Sigma(g)$
is $T_{g(p)}\R^{6k-1}/dg(T_pM^{4k})$. 
\begin{lma}\label{orxi}
The three line bundles $\det(T\tilde\Sigma(g))$, $\det(\xi(g))$, and
$\ker(dg)$ over $\tilde\Sigma(g)$ are isomorphic. The total space
$E(g)$ of the vector bundle 
$\xi(g)$ is orientable. Moreover, orientations on $M^{4k}$ and
$\R^{6k-1}$ induce an orientation on $E(g)$.
\end{lma}    
\begin{pf}
If $p$ is a point in $\tilde\Sigma(g)$ then,
$$
T_{g(p)}\R^{6k-1}=dg(T_p M^{4k}/\ker(dg)_p)\oplus \xi(g)_p.
$$
Hence, over $\tilde\Sigma(g)$, we have the following isomorphism of
line bundles:  
$$
\det(g^{\ast}(T\R^{6k-1}))=\det(TM^{4k})\otimes\ker(dg)
\otimes\det(\xi(g)).
$$
It follows that orientations on $M^{4k}$ and $\R^{6k}$ induce an
isomorphism 
\begin{equation}\label{is1}
\det(\xi(g))\approx\ker(dg).
\end{equation}

Consider a point $p\in\tilde\Sigma(g)$ which is not a cusp. Let $U$ be
a small neighborhood of $p$ in $M^{4k}$ in which there are no cusps. 
Let $g'=g|U$. If $U$ is chosen small enough then $g'$ does not have any
triple points, $\tilde\Sigma(g')=\tilde\Sigma(g)\cap U$, 
$X(g')=\tilde D_2(g')\cup \tilde\Sigma(g')$ is an embedded
submanifold of $U$, and $\tilde\Sigma(g')$ divides $X(g')$ into two
connected components $\tilde D_2^+(g')$ and $\tilde D_2^-(g')$. 

Let $x\in \tilde D_2(g')$. Then there exists a
unique $y\in\tilde D_2(g')$ such that $g'(x)=g'(y)$. Moreover,
there are neighborhoods $V$ of $x$ and $W$ of $y$ in $U$ such that
$A(x)=g'(V)\cap g'(W)\subset\R^{6k-1}$ is canonically diffeomorphic to
a neighborhood $B(x)$ of $x$ in $X(g')$. We pull back the
orientation induced on $A(x)$ by considering it as the {\em ordered}
intersection $g'(V)\cap g'(W)$ to $B(x)$. In this way we get an
orientation of $X(g')-\tilde\Sigma(g)$.

Note that the orientation on $A=A(x)=A(y)$ pulled back to the
neighborhood $B(x)\subset X(g')$ of $x$ is  
{\em opposite} to that pulled back to the neighborhood 
$B(y)\subset X(g')$ of $y$. Thus, the orientations induced on
$\tilde\Sigma(g')$ as the boundary of $ \tilde D_2^+(g')$ and 
$\tilde D_2^-(g')$,
respectively are {\em opposite}. This means that the orientation of
$\tilde D_2(g')$ extends over $\tilde\Sigma(g')$ and gives an
orientation of $X(g')$. 

Along $\tilde\Sigma(g')$, we have 
$T X(g')|\tilde\Sigma(g')= T\tilde\Sigma(g')\oplus\ker(dg')$. 
This is a subbundle of $T M^{4k}|\tilde\Sigma(g')$. 
Let $l$ be a $(2k-1)$-dimensional subbundle of 
$T M^{4k}|\tilde\Sigma(g')$ complementary to 
$T X(g')|\tilde\Sigma(g')$, i.e. such that, there is a bundle 
isomorphism  
$$
T M^{4k}|\tilde\Sigma(g')=T X(g')|\tilde\Sigma(g')\oplus l.
$$ 
Then orientations of $X(g')$ and $T M^{4k}$ induce an orientation on $l$.
The decomposition 
$$
T M^{4k}|\tilde\Sigma(g')=T\tilde\Sigma(g')\oplus\ker(dg')\oplus l,
$$    
gives the line bundle isomorphism
$$
\det(T M^{4k}|\tilde\Sigma(g'))=
\det(T\tilde\Sigma(g'))\otimes\ker(dg')\otimes\det(l).
$$
The orientations on $M^{4k}$ and $l$ then induce an isomorphism
$$
\det(T\tilde\Sigma(g'))\approx\ker(dg').
$$
Now $T_p\tilde\Sigma(g')=T_p\tilde\Sigma(g)$ and 
$\ker(dg')_p=\ker(dg)_p$. Hence, we have an isomorphism
\begin{equation}\label{is2}
\det(T\tilde\Sigma(g))\approx\ker(dg),
\end{equation}
over $\tilde\Sigma(g)-C$, where $C$ denotes the set of cusp points. Since
$C$ has codimension $2k$ in $\tilde\Sigma(g)$ the isomorphism extends
uniquely over $C$.

Finally, if $v$ is a point in $E(g)$, with $\pi(v)=p$, where $\pi\colon
E(g)\to\tilde\Sigma(g)$ is the projection. Then the tangent space of
$E(g)$ at $v$ splits as 
$$
T_vE(g)\approx T_p\tilde\Sigma(g)\oplus\xi(g)_p.
$$
Hence,
$$
\det(TE(g))\approx\ker(dg)\otimes\ker(dg), 
$$
is trivial and $E(g)$ is orientable. Moreover, the isomorphisms
~\eqref{is1} and ~\eqref{is2} give a trivialization of 
$\det(TE(g))$. That is, there is an induced orientation on $E(g)$.
\end{pf}
\begin{rmk}
Alternatively, one can prove Lemma ~\ref{orxi} by
using the diagram
$$
\begin{CD}
M^{4k}\times M^{4k}  @>>{g\times g}> \R^{6k-1}\times\R^{6k-1}\\
@V{\phi}VV  @VV{\phi}V \\
M^{4k}\times M^{4k}  @>>{g\times g}> \R^{6k-1}\times\R^{6k-1}
\end{CD},
$$
where $\phi(x,y)=(y,x)$, to show that the closure $X(g)$ of the double
point manifold of $g$ in $M^{4k}$ is oriented.   
With this shown, one can use the identification of
the normal bundle of $\tilde\Sigma(g)-C$ ($C$ is the set of cusp
points) in $X(g)$ with $\ker(dg)$ to establish the isomorphism 
$\det(T\tilde\Sigma(g))\approx \ker(dg)$.
\end{rmk}

Lemma ~\ref{orxi} says that $\xi(g)$ is a $2k$-dimensional vector
bundle over the $2k$-dimensional manifold $\tilde\Sigma(g)$ with total
space $E(g)$ which is oriented. There is an Euler number associated to
such a bundle. In terms of homology with local coefficients this Euler
number is derived as follows: 

Let ${\mathcal K}$ be the twisted integer
local coefficient system associated to the orientation bundle of
$\tilde\Sigma$ and let ${\mathcal F}$ denote the local system
associated to the fiber orientation of the bundle $\xi(g)$. We have
proved above that 
${\mathcal K}\approx{\mathcal F}$ (both are isomorphic to the sheaf
of unit length sections of $\ker(dg)$). The orientation on the total
space gives the relation 
\begin{equation}\label{eqtwist}
{\mathcal K}\otimes{\mathcal F}=\Z,
\end{equation}     
where $\Z$ is the (trivialized) local coefficient system associated to
the orientation bundle of $E(g)$ restricted to the zero-section. This
relation is used to specify an isomorphism 
${\mathcal K}\approx{\mathcal F}$, by requiring that at each point it
is given by ordinary multiplication of integers. This specified
isomorphism gives in turn a well-defined pairing
$$
H^{2k}(\tilde\Sigma(g);{\mathcal F})\otimes
H_{2k}(\tilde\Sigma(g),{\mathcal K})\to\Z.
$$
Especially, we get the Euler number 
$e(\xi(g))=\langle e,[\tilde\Sigma(g)]\rangle$, where 
$e\in H^{2k}(\tilde\Sigma(g);{\mathcal F})$ is the Euler class of
$\xi(g)$ and 
$[\tilde\Sigma(g)]\in H_{2k}(\tilde\Sigma(g),{\mathcal K})$ is the 
orientation class of $\tilde\Sigma(g)$.

One can compute the above Euler number by choosing a section $s$ of
$\xi(g)$ transverse to the zero-section and sum up the local
intersection numbers at the zeros of $s$. The local intersection number
at a zero $p$ of $s$ is the intersection number in $E(g)$ of the
zero-section with some local orientation on a neighborhood $U$ of $p$
and $s(U)$ with the orientation induced from that chosen local
orientation on $U$.

\subsection{Prim maps}
\begin{dfn}
Let $M^n$ be a manifold. A map $g\colon M^n\to\R^{n+m}$ is a {\em prim
map (projected immersion)} if there exists an immersion $G\colon
M^{n}\to\R^{n+m+1}$ such 
that $g=\pi\circ G$, where $\pi\colon\R^{n+m+1}\to\R^{n+m}$ is the
projection forgetting the last coordinate.
\end{dfn}

\begin{rmk}\label{xiprim}
The bundle $\xi=\xi(g)$ has a more simple description if $g\colon
M^{4k}\to\R^{6k-1}$ is a prim map. 
Assume that $g=\pi\circ G$. Then $\xi$ is isomorphic to
the normal bundle $\nu$ of $G$ restricted to $\tilde\Sigma(g)$. 

Since the homology class of $\tilde\Sigma(g)$ in $M^{4k}$ is dual to
the normal 
Euler class of the immersion $G$, we have (if $e(\eta)$ denotes the
Euler class of the bundle $\eta$)  
$$
\langle e(\xi),[\tilde\Sigma]\rangle=
\langle e^2(\nu),[M^{4k}]\rangle=\bar p_k[M^{4k}].
$$
\end{rmk}
\subsection{The Euler class and cusps}
In this section we prove that the expressions in Equations ~\eqref{grom1}
and ~\eqref{grom2} agree. We use the notation of Section ~\ref{not1}.
\begin{lma}\label{s11clo}
If $g\colon M^{4k}\to\R^{6k-1}$ is a generic map of a closed oriented
manifold then
$$
\bar p_k(M^{4k})=\sharp\Sigma^{1,1}(g).
$$
\end{lma}  
\begin{pf}
This is Lemma 3 in ~\cite{S2}. 
\end{pf}

\begin{lma}\label{es11}
Let $M^{4k}$ be a closed oriented manifold and let 
$g\colon M^{4k}\to\R^{6k-1}$ be a generic map. Then 
\begin{equation}\label{es11eq}
e(\xi(g))=\sharp\Sigma^{1,1}(g).
\end{equation}
\end{lma}
\begin{pf}
If $g$ is a prim map then both sides in Equation ~\eqref{es11eq} equal
$\bar p_k[M^{4k}]$, see Remark ~\ref{xiprim} for the left-hand side
and Lemma ~\ref{s11clo} for the right-hand side. 
Both sides are invariant under cobordisms of generic maps and hence they
define homomorphisms $\Omega_{4k}(\R^{6k-1})\to\Z$. By Burlet's
theorem ~\cite{Bu} (see the proof of Lemma ~\ref{lma1}) the classes
representable by prim maps form a subgroup of finite index. Since the
homomorphisms agree on this subgroup they agree on the whole group.
\end{pf}

\begin{rmk}\label{newsgn}
Alternatively, Lemma ~\ref{es11} may be proved as follows:

The second derivative $D^2g$ (intrinsic derivative of Porteous) of the
map $g$ is a quadratic form on 
$\ker(dg)$ with values in $\xi(g)=T\R^{6k-1}/\img(dg)$. Choose a Riemannian
metric on $\ker(dg)$ and let $\pm v(x)$ denote the unit vectors in
each fiber over $x\in\Sigma(g)$. Then $s(x)=D^2g(\pm v(x),\pm v(x))$ gives
a section of $\xi(g)$. This section vanishes exactly at the
$\Sigma^{1,1}$-points. (Confer Section 2.2 of ~\cite{Hae2} which
deals with the cusp-free case.)  
The Euler number of the bundle $\xi(g)$ is the sum of
local intersection numbers at the zeros of $s$. 
\end{rmk}

\subsection{Cobordism groups of maps and natural homomorphisms}
We shall consider classifying spaces of certain maps of
codimension $m$ into Euclidean space.
\begin{dfn}\label{clsp}
\begin{itemize}
\item[{\rm (a)}]
Let $X(m)$ be the classifying space of (generic) codimension $m$ maps
of closed oriented manifolds into Euclidean space.
\item[{\rm (b)}]
Let $\bar X(m)$ be the classifying space of codimension $m$ prim maps
of closed oriented manifolds into Euclidean space.
\item[{\rm (c)}]
Let $\Gamma(m)$ be the classifying space of codimension $m$ immersions
of closed oriented manifolds into Euclidean space.
\item[{\rm (d)}]
Let $\Gamma_{{\rm fr}}(m)$ be the classifying space of codimension $m$
framed immersions of closed oriented manifolds into Euclidean space.
\end{itemize}
\end{dfn}
Definition ~\ref{clsp} (a) means that $\pi_{n+m}(X(m))$
is the cobordism group of arbitrary 
maps of oriented $n$-manifolds in $\R^{n+m}$. Definition
~\ref{clsp} (b)-(d) give similar interpretations of the homotopy
groups of the corresponding spaces. 

\begin{rmk}
Up to homotopy equivalence these spaces can be identified as follows:
\begin{itemize}
\item[{\rm (a)}]
$$
X(m)=\lim_{K\to\infty}\Omega^{K}MSO(K+m),
$$ 
\item[{\rm (b)}]
$$
\bar X(m)=\lim_{K\to\infty}\Omega^{K+1}\Su^{K}MSO(m+1),
$$ 
\item[{\rm (c)}]
$$
\Gamma(m)=\lim_{K\to\infty}\Omega^{K}\Su^{K}MSO(m),
$$ 
\item[{\rm (d)}]
$$
\Gamma_{{\rm fr}}(m)=\lim_{K\to\infty}\Omega^{K}\Su^{K+m},
$$ 
\end{itemize}
where $\Omega^{j}$ denotes the $j^{\rm th}$ loop space and
$\Su^{j}$ the $j^{\rm th}$ suspension or the $j$-dimensional sphere.
\end{rmk}

We note that there are natural inclusions among these spaces and that
the corresponding relative homotopy groups also have concrete geometric
interpretations. For example, 
$\pi_{n+m}(\bar X(m),\Gamma_{{\rm fr}}(m))$ is the cobordism group of
prim maps $(M^n,\partial M^n)\to(\R^{n+m}_+,\partial\R^{n+m}_+)$ which
are framed immersions on the boundary.  

\begin{dfn}\label{dfnhomS11}
Let $\beta$ be an element of $\pi_{6k-1}(\bar X(2k-1))$ or of 
$\pi_{6k-1}(\bar X(2k-1),\Gamma_{{\rm fr}}(2k-1))$, or of
$\pi_{6k-1}(X(2k-1),\Gamma_{{\rm fr}}(2k-1))$. Represent $\beta$ by a
generic map $g\colon M^{4k}\to\R^{6k-1}$ and define 
$$
\Sigma^{1,1}(\beta)=\sharp\Sigma^{1,1}(g),
$$ 
see Definition ~\ref{dfnS11}. 
\end{dfn}
Note that a generic cobordism between two generic maps as in Definition
~\ref{dfnhomS11} gives a cobordism between their $0$-dimensional
manifolds of cusp points. Hence, the function $\Sigma^{1,1}$ is a well
defined homomorphism. 

\begin{dfn}\label{dfnhom3t}
Let $\beta$ be an element of $\pi_{6k}(\Gamma(2k))$.
Represent $\beta$ by a
generic (self transverse) immersion $g\colon M^{4k}\to\R^{6k}$ and define 
$$
3t(\beta)=3 t(g),
$$ 
see Definition ~\ref{dfnt}. 
\end{dfn}
As above we see that $3t$ is a well defined homomorphism.

\begin{dfn}\label{dfnhom3t+L}
Let $\beta$ be an element of $\pi_{6k}(\Gamma(2k),\Gamma_{{\rm
fr}}(2k))$. 
Represent $\beta$ by a
generic (self transverse) immersion $g\colon (M^{4k},\partial
M^{4k})\to(\R^{6k}_+,\partial\R^{6k}_+)$ and define 
$$
\Phi(\beta)=3 t(g)+L(\partial g),
$$ 
see Definitions ~\ref{dfnt} and ~\ref{dfnL}. 
\end{dfn}
\begin{dfn}\label{dfnhom3t+L-l}
Let $\beta$ be an element of $\pi_{6k}(X(2k),\Gamma_{{\rm fr}}(2k))$.
Represent $\beta$ by a generic map 
$g\colon (M^{4k},\partial M^{4k})\to(\R^{6k}_+,\partial\R^{6k}_+)$ and
define  
$$
\Psi(\beta)=3 t(g)-3 l(g) +L(\partial g),
$$ 
see Definitions ~\ref{dfnt}, ~\ref{dfnL}, and ~\ref{dfnl}. 
\end{dfn}
By the same argument as is used in the proof of Theorem ~\ref{thm1}
(b) to show that $\Theta$ (see Equation ~\eqref{wdexpr}) is well
defined, it follows that $\Phi$ and $\Psi$ as defined above are well
defined homomorphisms.  

The following lemma is the main step in the proof of Theorem
~\ref{thm1} (a).
\begin{lma}\label{diagram}
Let 
$$
i\colon \pi_{6k-1}(X(2k-1),\Gamma_{{\rm fr}}(2k-1))\to
\pi_{6k}(X(2k),\Gamma_{{\rm fr}}(2k))
$$
be the natural homomorphism. Then $\Sigma^{1,1}=\Psi\circ i$.
\end{lma}
\begin{pf}
To shorten the notation, let $m=2k-1$ and $n=2k$. Consider the
following diagram:
$$
\begin{CD}
\pi_{6k-1}(\bar X(m)) @>>> 
\pi_{6k-1}(\bar X(m),\Gamma_{{\rm fr}}(m)) @>>> 
\pi_{6k-1}(X(m),\Gamma_{{\rm fr}}(m))\\ 
@V{i''}VV @V{i'}VV  @VViV\\
\pi_{6k}(\Gamma(n)) @>>> 
\pi_{6k}(\Gamma(n),\Gamma_{{\rm fr}}(n)) @>>> 
\pi_{6k}(X(n),\Gamma_{{\rm fr}}(n))
\end{CD}
$$

The horizontal homomorphisms are obtained by forgetting 
structure. The vertical homomorphisms are the natural ones induced by
the (framed) inclusion $\R^{6k-1}\to\R^{6k}$. The diagram clearly
commutes. 

We now show that the horizontal homomorphisms in these sequences all
have finite cokernels. We start with the first horizontal arrow. This
homomorphism is part of the long exact homotopy sequence of the pair
$(\bar X(m),\Gamma_{{\rm fr}}(m))$. Consider the following fragment of this
sequence:  
$$
\begin{CD}
\dots \pi_{6k-1}(\bar X(m)) @>>> 
\pi_{6k-1}(\bar X(m),\Gamma_{{\rm fr}}(m)) @>>> 
\pi_{6k-2}(\Gamma_{{\rm fr}}(m)) \dots
\end{CD}
$$
The group $\pi_{6k-2}(\Gamma_{{\rm fr}}(m))$ is the stable
homotopy group of spheres $\pi^{\rm s}(4k-1)$, which is finite. It
follows that the cokernel is finite.

Similarly, if the  groups
\begin{itemize}
\item[(a)]$\pi_{6k-1}(X(m),\bar X(m))$,
\item[(b)]$\pi_{6k-1}(\Gamma_{{\rm fr}}(n))$, and
\item[(c)]$\pi_{6k}(X(n),\Gamma(n))$,
\end{itemize}
are finite then the other horizontal homomorphisms have finite
cokernels. 

The group in (b) is again a stable homotopy group of
spheres and therefore finite. The group in (c) is finite by Burlet's
theorem ~\cite{Bu} (see the proof of Lemma ~\ref{lma1}). 
Finally, the group in (a) is isomorphic to the group in (c): 

Any map 
$f\colon (M^{4k},\partial M^{4k})\to(\R^{6k-1}_+,\partial\R^{6k-1}_+)$
which is prim on the boundary lifts to a map 
$F\colon (M^{4k},\partial M^{4k})\to(\R^{6k}_+,\partial\R^{6k}_+)$
which is an immersion on the boundary. This gives the isomorphism 
$\pi_{6k-1}(X(m),\bar X(m))\to\pi_{6k}(X(n),\Gamma(n))$ on
representatives. The inverse is induced by the projection 
$(\R^{6k}_+,\partial\R^{6k}_+)\to(\R^{6k-1}_+,\partial\R^{6k-1}_+)$.

By Remark 1 in ~\cite{S2}, it follows that 
$\Sigma^{1,1}=3t\circ i''$. Since 
all cokernels are finite, we conclude first that
$\Sigma^{1,1}=\Phi\circ i'$ 
and then that $\Sigma^{1,1}=\Psi\circ i$. 
\end{pf}  
\subsection{Proof of Theorem ~\ref{thm1} (a)}\label{pfthm1a}
Let $f\colon S^{4k-1}\to\R^{4k+1}$ be an immersion. Then there is a
homotopically unique normal framing of $f$. Thus, the immersion 
$j\circ f\colon S^{4k-1}\to\R^{6k-1}$ is also framed. 

Let $g\colon M^{4k}\to\R^{6k-1}$ be a generic map of a manifold with
spherical boundary such that $\partial g$ is regularly homotopic to
$j\circ f$. Then there is an induced normal framing of $\partial g$. 
After composing $g$ with a translation, we can assume that
$g(M^{4k})$ is contained in the half space on which the last
coordinate function is strictly positive.

Applying Hirsch lemma to any vector field in the normal framing of
$\partial g$, we find a homotopy of $M^{4k}$ supported in a small
collar of the boundary 
$\partial M^{4k}$, which is a regular homotopy when restricted to this
collar 
and a framed regular homotopy of $\partial M^{4k}$ and which deforms
$\partial g$ to an immersion mapping into $\partial \R^{6k-1}_+$. 

Let $g'\colon (M^{4k},\partial
M^{4k})\to(\R^{6k-1}_+,\partial\R^{6k-1}_+)$ denote the map obtained
from $g$. Then $g'$ represents an element
$\zeta\in\pi_{6k-1}(X(2k-1),\Gamma_{{\rm fr}}(2k-1))$. 

Theorem ~\ref{thm1} (b) says that
$a_k(2k-1)!\cdot\Omega(f)+\bar p_k[\hat M^{4k}]=\Psi(i(\zeta))$. Thus,
by Lemma ~\ref{diagram}, we have
$$
a_k(2k-1)!\cdot \Omega(f)+\bar p_k[\hat M^{4k}]=
\Sigma^{1,1}(\zeta)=\sharp\Sigma^{1,1}(g')=\sharp\Sigma^{1,1}(g).
$$

This proves Equation ~\eqref{grom1}. Equation ~\eqref{grom2} then follows
from Lemma ~\ref{es11}.\qed

\section{$\Imm(S^{4k-1},\R^{4k+1})\to\Imm(S^{4k-1},\R^{6k-1})$ is
injective}
In this section we prove Corollary ~\ref{corthm1}.
\subsection{Proof of Corollary ~\ref{corthm1}}\label{pfcorthm1}
Let $f,g\colon S^{4k-1}\to\R^{4k+1}$ be immersions such that 
$j\circ f\colon S^{4k-1}\to\R^{6k-1}$ is regularly homotopic to
$j\circ g$. Theorem ~\ref{thm1}
applied to a 
generic immersion $\partial h\colon S^{4k-1}\to\R^{6k-1}$, in the
common regular
homotopy class of $j\circ f$ and $j\circ g$, bounding a generic
map $h\colon M^{4k}\to \R^{6k}_+$, shows that
$\Omega(f)=\Omega(g)$. In other words, 
$\Imm(S^{4k-1},\R^{4k+1})\to\Imm(S^{4k-1},\R^{6k-1})$ is injective.

Let $f,g\colon S^{4k-1}\to\R^{4k+1}$ be immersions such that
$\Omega(f)-\Omega(g)=c$. Let $h\colon M^{4k}\to\R^{6k-1}$ be a
generic map such that $\partial h=j\circ g$. If 
$F\colon S^{4k-1}\times I\to\R^{6k-1}$ is a generic homotopy connecting
$j\circ f$ to $j\circ g$ then we glue the collar 
$S^{4k-1}\times I$ to $M^{4k}$ and obtain a manifold $N^{4k}$ with a generic
map $k\colon N^{4k}\to\R^{6k-1}$, which equals $F$ on the collar and $h$
on $M^{4k}$. Now, $\partial k=j\circ f$ and $\hat M^{4k}\approx \hat
N^{4k}$. Hence, by Theorem ~\ref{thm1} (a)
$$
c\cdot a_k(2k-1)!=(\Omega(f)-\Omega(g))\cdot a_k(2k-1)!
=\sharp\Sigma^{1,1}(k)-\sharp\Sigma^{1,1}(h)=\sharp\Sigma^{1,1}(F).
$$
\qed

\section{On the number of $(2k+1)$-tuple points of an immersion
$M^{4k} \to R^{4k+2}_+$ bounding a given immersion}\label{mpscob} 
In this section we prove Theorem ~\ref{thm3}.
\subsection{Proof of Theorem ~\ref{thm3}}\label{pfthm3}
Let $\Theta(f)$ denote the left-hand side of Equation ~\eqref{high0}. We
shall show that
\begin{itemize}
\item[(a)] $\Theta(f)$ is well-defined, i.e. it does not depend on the
choices of $d$ and $g$. 
\item[(b)] $\Theta(f)$ does not change if the immersion $f$ changes in
its cobordism class.
\end{itemize}

But first we show that the theorem follows from (a) and (b):

Let $\Imm^{SO}(4k-1,2)$ be the cobordism group of immersions of oriented
$(4k-1)$-manifolds in $\R^{4k+3}$. Then $\Imm^{SO}(4k-1,2)$ is
isomorphic to $\pi^{\rm s}_{4k+3}(\C P^{\infty})$, the $(4k+3)^{{\rm th}}$
stable homotopy group of $\C P^{\infty}$. By Serre's theorem, 
$\pi_n^{\rm s}(X)\otimes\Q \approx H_n(X)\otimes\Q$, and hence,
$\Imm^{SO}(4k-1,2)$ is finite.    

Let $\mathcal I_2$ be the subgroup of $\Imm^{SO}(4k-1, 2)$ which consists
of cobordism classes representable by immersions of $2$-connected 
manifolds. By (a) and (b), the formula $\Theta(f)$ defines a
homomorphism $\Theta\colon\mathcal I_2\to\Z$. Since $\mathcal I_2$ is
finite, this homomorphism must be identically zero.

We now return to the proofs of (a) and (b).

The statement (a) follows from a special case of Herbert's theorem
~\cite{He} 
saying that if $h\colon N^{4k}\to\R^{4k+2}$ is an immersion of a
{\em closed} manifold then $\la\bar p_1^k,[N^{4k}]\ra=\sharp D_{2k+1}(h)$. 

Consider (b). Let $F\colon W^{4k}\to\R^{4k+1}\times I$ be a generic
immersion, which is a cobordism between 
$f_0\colon V_0\to\R^{4k+1}\times 0$ and 
$f_1\colon V_1\to\R^{4k+1} \times 1$. We must show.
$$
L_{2k}(f_1) - L_{2k}(f_0) = \bar p_1^k[W^{4k}] + 
(2k+1)\sharp D_{2k+1}(F). 
$$

We first introduce some notation: Let $\Delta_{2k}(F)$ denote the
resolved $2k$-fold self intersection manifold of $F$. (Then there is a
map $\Delta_{2k}(F)\to D_{2k}(F)\cup D_{2k+1}(F)$, which is a
diffeomorphism when restricted to the preimage of $D_{2k}(F)$ and in
the preimage of each point in $D_{2k+1}(F)$ there are exactly $2k+1$
points.) Similarly, let $\tilde\Delta_{2k}(F)$ denote the resolution
of $\tilde D_{2k}(F)\cup \tilde D_{2k+1}(F)$. Then there is a 
$2k$-fold covering
$\pi\colon\tilde\Delta_{2k}(F)\to\Delta_{2k}(F)$. Also, let
$j\colon\tilde\Delta_{2k}(F)\to W^{4k}$ and
$i\colon\Delta_{2k}(F)\to\R^{4k+2}\times I$ denote the obvious
immersions.  
 
Let $\nu_F$ denote the normal bundle of $F$ and let
$\zeta=j^*(\nu_F)$. The $2$-connectedness of $V_0$ and $V_1$ implies 
that $\nu_F|\partial W^{4k}$ has a homotopically unique
trivialization. The same is then true for the restriction of $\zeta$
to the boundary. Let $s$ be a section of $\zeta$ which does not vanish
on the boundary. 

The normal bundle of $i\colon\Delta_{2k}(F)\to\R^{4k+1}\times I$ is then
$\pi_{!}(\zeta)$ and the section $s$ gives a section $z$ of
$\pi_{!}(\zeta)$, namely $z(q)=s(q_1)+\dots+s(q_{2k})$, where
$q=\pi(q_1)=\dots=\pi(q_{2k})$.  

Let $i'\colon\Delta_{2k}(F)\to\R^{4k+2}\times I$ be given by
$i$ shifted a small distance along
$z$ (i.e. $i'(x)=i(x)+\epsilon z(x)$, where $\epsilon>0$ is very small). 
Then $i'(\Delta_{2k}(F))\cap F(W^{4k})$ is a collection of points. The
points are of two types:
\begin{itemize}
\item[(i)]
Near each $(2k+1)$-fold self intersection point $p$ of $F$ there are
$2k+1$ intersection points, all with the same local intersection
number, which is the sign of the $(2k+1)$-fold self intersection point
$p$.
\item[(ii)]
One intersection point for each zero of $s$. There is a local
intersection number associated to such a point.  
\end{itemize}
Clearly, we can choose $s$ so that the sets of intersection points of
type (i) and (ii), respectively, are disjoint.

Hence, 
$$
i'(\Delta_{2k}(F))\bullet F(W^{4k})=(2k+1)\sharp D_{2k+1}(F)+\sharp\{
s^{-1}(0)\},
$$ 
where $\sharp\{s^{-1}(0)\}$ denotes the algebraic number of
intersection points of type (ii). Applying Lemma ~\ref{lkint}, we find that
$$
L_{2k}(f_1)-L_{2k}(f_0)=(2k+1)\sharp D_{2k+1}(F)+\sharp\{ s^{-1}(0)\}.
$$  
Now, $\sharp \{s^{-1}(0)\}$ equals the relative Euler class 
$e(\zeta)\in H^2(\tilde\Delta_{2k}(F),\partial\tilde\Delta_{2k}(F))$
of the bundle $\zeta$. The immersed 
submanifold  $j(\tilde\Delta_{2k}(F))$ represents
the class dual to $(2k-1)^{\rm th}$ power of the Euler class of $\nu(F)$.
If $\mathcal D$ is the Poincar\'e duality operator on $W=W^{4k}$ then
\begin{align*}
\sharp \{s^{-1}(0)\} &=
\left\langle e(j^\ast\nu_F), [\tilde \Delta_{2k}(F), \partial \tilde
\Delta_{2k}(F)]\right\rangle =
\left\langle j^\ast e(\nu_F), [\tilde \Delta_{2k}(F), \partial \tilde
\Delta_{2k}(F)]\right\rangle =\\
&= \left\langle e(\nu_F), j_\ast[\tilde \Delta_{2k}(F), \partial \tilde
\Delta_{2k}(F)]\right\rangle = 
\left\langle e(\nu_F), {\mathcal D} e^{2k-1}(\nu(F)) \right\rangle =\\
&= \left\langle e^{2k}(\nu_F), [W, \partial W]\right\rangle
= \langle\bar p_1^k,[W,\partial W]\rangle.
\end{align*}
\qed

\section{Codimension two immersions of spheres\\ of dimensions $8k+5$
and $8k+1$}\label{sist}
In this section we state a formula which might give the Smale
invariant of an 
immersion $S^{8k+1}\to\R^{8k+3}$ and prove that the corresponding
formula vanishes identically 
for immersions $S^{8k+5}\to\R^{8k+7}$.
\subsection{A brief discussion of definitions and notation}\label{dfnnot}
Let $f\colon S^{8k+1}\to\R^{8k+3}$ ($f\colon S^{8k+5}\to\R^{8k+7}$) be
an immersion. Let 
$j\colon\R^{8k+3}\to\R^{12k+2}$ ($j\colon\R^{8k+7}\to\R^{12k+8}$) be
the inclusion. Let
$g\colon M^{8k+2}\to\R^{12k+3}_+$ ($g\colon M^{8k+6}\to\R^{12k+9}_+$)
be a generic map of a compact 
manifold such that $\partial g$ is an immersion regularly homotopic to
$j\circ f$.

The {\em $\Z_2$-valued} invariants $t(g)$, $l(g)$, and  
$L(\partial g)$ are then defined as the corresponding invariants in
Definitions ~\ref{dfnt}, ~\ref{dfnl}, and ~\ref{dfnL}, respectively,
with the following modifications: First,
$\Z_2$ is used instead of $\Z$ and second, no orientations are needed.

For a compact closed $2j$-dimensional manifold $M^{2j}$, let  
$\bar w_j$ denote the $j^{\rm th}$ normal Stiefel-Whitney class
of $M^{2j}$ and let $\bar w_j^2[M^{2j}]\in\Z_2$ denote the
corresponding Stiefel-Whitney number.
\subsection{Possibly a Smale invariant formula}
Let $f\colon S^{8k+1}\to\R^{8k+3}$ be an immersion. Let
$\Omega(f)\in\Z_2$ be its Smale invariant and let $g$ be as in Section
~\ref{dfnnot}.
\begin{cnj}\label{cnj1}
Either
\begin{equation}\label{cnjfmla}
\Omega(f)=w_{4k+1}^2[\hat M^{8k+2}]+t(g)+l(g)+L(\partial g),
\end{equation}
or
\begin{equation}\label{nullfmla}
w_{4k+1}^2[\hat M^{8k+2}]+t(g)+l(g)+L(\partial g)=0,
\end{equation}
where we use the notations introduced in Section ~\ref{dfnnot}. Which
one is true? (Note that Equations ~\eqref{cnjfmla} and
~\eqref{nullfmla} are equations in $\Z_2$.)
\end{cnj}

Question ~\ref{cnj1} can be treated in the same way as Theorem
~\ref{thm1} (b). In the proof of Theorem ~\ref{thm1} (b) we used Lemma
~\ref{lma1}. The analog of Lemma ~\ref{lma1} in the present
situation is the following:

\begin{lma}\label{mod2}
Let $M^{2n}$ be a closed $2n$-dimensional manifold and let 
$h\colon M^{2n}\to\R^{3n}$ be any generic map. Then
$w_{n}^2(M^{2n})+t(h)+l(h)=0$ (in $\Z_2$).
\end{lma}
\begin{pf}
Lemma ~\ref{mod2} is proved  in ~\cite{S}.
\end{pf}

It follows from Lemma ~\ref{mod2} and the fact that $L(\partial g)$
changes at triple point instances of generic regular homotopies (see
~\cite{E3}), that the right-hand side of
~\eqref{cnjfmla} is independent of $g$. 
Hence, the right-hand side of ~\eqref{cnjfmla} induces a homomorphism
$\Imm(S^{8k+1},\R^{8k+3})\to\Z_2$ which is either zero (Equation
~\eqref{nullfmla} is true) or the Smale
invariant (Equation ~\eqref{cnjfmla} is true). 

Thus, to prove ~\eqref{cnjfmla} it is enough to find
one example for which the right-hand side of ~\eqref{cnjfmla} does not
vanish. 

\begin{rmk}
In the same way that Corollary ~\ref{corthm1} follows from Theorem
~\ref{thm1}, it would follow from Equation ~\ref{cnjfmla} that 
$\Imm(S^{8k+1},\R^{8k+3})\to\Imm(S^{8k+1},\R^{12k+2})$ is
injective. For $k=1$ this might be the case. Indeed,
$\Imm(S^9,\R^{11})=\Z_2$ and $\Imm(S^9,\R^{14})=\Z_2$, see ~\cite{P}. 
\end{rmk}
\subsection{A formula expressing Smale invariant equals zero}
There is only one regular homotopy class of immersions
$S^{8k+5}\to\R^{8k+7}$. Hence, Lemma ~\ref{mod2} and the fact that $L$
changes under triple point instances of generic regular homotopies
imply the following (with notation as in Section ~\ref{dfnnot}):
\begin{prp}\label{prpsist}
Let $f\colon S^{8k+5}\to\R^{8k+7}$ be an immersion. Then
$$
w_{4k+3}^2[\hat M^{8k+6}]+t(g)+l(g)+L(\partial g)=0,
$$
where we use the notation introduced in Section ~\ref{dfnnot}. (Note
that this is an equation in $\Z_2$.)
\qed
\end{prp}

\end{document}